\documentclass[12pt,reqno]{amsart}

\usepackage{amsmath}

\usepackage{amsfonts}

\usepackage{hyperref}

\usepackage{graphicx}

\theoremstyle{plain}

\newtheorem{defi}{Definition}

\newcommand{\brdef}{\begin{defi}}

\newcommand{\erdef}{\end{defi}}

\newcommand{\brexam}{\begin{examp}}

\newcommand{\erexam}{\end{examp}}

\newtheorem{cor}{Corollary}

\newcommand{\bcor}{\begin{cor}}

\newcommand{\ecor}{\end{cor}}

\newtheorem{thm}{Theorem}

\newcommand{\bth}{\begin{thm}}

\newcommand{\eth}{\end{thm}}

\newtheorem{lem}{Lemma}

\newcommand{\ble}{\begin{lem}}

\newcommand{\ele}{\end{lem}}

\newtheorem{pro}{Proposition}

\newcommand{\bpro}{\begin{pro}}

\newcommand{\epro}{\end{pro}}

\huge

\numberwithin{equation}{section}

\begin{document}
\title{$\mathcal{Z}$-symmetries of $(\epsilon)$-para-Sasakian 3-manifolds}
\author{D. G. Prakasha$^{1, *}$, P. Veeresha$^{2}$ and M. Nagaraja$^{3}$}
\maketitle

\begin{center}
\pagestyle{myheadings} \markboth{D. G. Prakasha, P. Veeresha and Inan Unal}{Invariant Submanifolds of Generalized Sasakian-space-forms}
$^{1}$Department of Mathematics, Faculty of Science, Davangere University, Shivagangothri, Davangere - 577 007, India.\\[0pt]
e-mail: {\verb|prakashadg@gmail.com|}\\[0pt]
$^{2}$Department of Mathematics, Karnatak University, Dharwad - 580 003, India.\\[0pt]
e-mail: {\verb|prakashadg@gmail.com|, \verb|viru0913@gmail.com|}\\[0pt]
$^{3}$Department of Mathematics, Tunga Mahavidhyalaya, Thirthahalli - 577 432, India.\\[0pt]
e-mail: {\verb|nagarajtmvt@gmail.com|}
\end{center}
\begin{quotation}
\textbf{Abstract:}\,\  The object of this paper is study $(\epsilon)$-para-Sasakian 3-manifolds satisfying certain conditions on the $\mathcal{Z}$ tensor. We characterize, $\mathcal{Z}$-symmetric; $\mathcal{Z}$-semisymmetric; $\mathcal{Z}$-pseudosymmetric; and projectively $\mathcal{Z}$-semisymmetric conditions on an $(\epsilon)$-para-Sasakian 3-manifold.
\end{quotation}


\section{\bf{Introduction}}


In 1969, Takahashi \cite{TT} initiated the study of almost contact manifolds associated with an indefinite metrics. These indefinite almost contact metric manifolds are also called as $(\epsilon)$-almost contact metric manifolds. The study of indefinite metric manifolds is of interest from the standpoint of physics and relativity. Indefinite metric manifolds have been studied by several authors. In 1993, Bejancu and Duggal \cite{BD} introduced the concept of $(\epsilon)$-Sasakian manifolds. Some interesting properties of these manifolds was studied in the papers \cite{KLD}, \cite{DS}, \cite{KRN}, \cite{XX} and the references therein. In 2009, De and Sarkar \cite{DS1} introduced the concept of $(\epsilon)$-Kenmotsu manifolds and showed that the existence of new structure on an indefinite metrics influences the curvatures. \\ 
\indent Tripathi and his co-authors \cite{TKYK} initiated the study of $(\epsilon)$-almost para-contact metric manifolds, which is not necessarily Lorentzian. In particular, they studied $(\epsilon)$-para-Sasakian manifolds, with the structure vector field $\xi$ is spacelike or timelike according as $\epsilon = 1$ or $\epsilon = -1$. An $(\epsilon)$-almost contact metric manifold is always odd dimensional but an $(\epsilon)$-almost para-contact metric manifold could be even dimensional as well. Later, Perktas and his co-authors \cite{PKTK} studied $(\epsilon)$-para-Sasakian manifolds in dimension 3. 
\par In 2012, Mantica and Molinari \cite{MM} defined a generalized (0, 2) symmetric $\mathcal{Z}$ tensor given by
\begin{equation}\label{1.1}
\mathcal{Z}(X, Y) = S(X, Y) + \psi g(X, Y),
\end{equation}
where $\psi$ is an arbitrary scalar function. Properties of $\mathcal{Z}$ tensor were pointed out in the papers \cite{MS} and \cite{MS1}. This tensor is a general notion of the Einstein gravitational tensor in General Relativity. Recently, Mallick and De \cite{MD} studied various properties of $\mathcal{Z}$ tensor on N(k)-quasi-Einstein manifolds.
\par The present paper is organized as follows: After preliminaries, in Section 3, we study $\mathcal{Z}$-semisymmetric, $\mathcal{Z}$-pseudosymmetric and projectively $\mathcal{Z}$-semisymmetric $(\epsilon)$-para-Sasakian 3-manifolds. Here, we prove that, for an $(\epsilon)$-para-Sasakian 3-manifold, the conditions of being $\mathcal{Z}$-semisymmetric; Ricci-symmetric; Ricci-semi symmetric; or Einstein manifold are all equivalent. Also show that, if an $(\epsilon)$-para-Sasakian 3-manifold $M$ is $\mathcal{Z}$-pseudosymmetric then it is either Ricci-semisymmetric or pseudosymmetric function $L_{\mathcal{Z}} = -\epsilon$ holds on $M$. Further, we prove that a projectively $\mathcal{Z}$-semisymmetric $(\epsilon)$-para-Sasakian 3-manifold is an Einstein manifold. 
\section{\bf Preliminaries}
A manifold $M$ is to admit an almost para-contact structure if it admit a tensor  field $\phi$ of type (1,1), a vector field $\xi$ and a 1-form $\eta$ satisfying
\begin{equation}
\label{2.1}
\phi^2 = I - \eta \otimes \xi, \,\,\eta(\xi)=1,\,\ \phi \xi = 0, \,\, \eta \cdot \phi=0.
 \end{equation} 
Let $g$ be a semi-Riemannian metric with index$(g)=v$ such that
\begin{equation}
\label{2.2}
 g(\phi X, \phi Y) = g(X, Y) - \epsilon \eta(X)\eta(Y). \,\,\, X, Y \in TM,
 \end{equation} 
where, $\epsilon=\pm 1$. Then $M$ is called  an  $(\epsilon)$-almost para-contact metric manifold equipped with an {\it $(\epsilon)$-almost para-contact metric structure} $(\phi, \xi, \eta, g, \epsilon)$. In particular, if index of $g$ is equal to one, then an $(\epsilon)$-almost para-contact metric manifold is said to be a {\it Lorentzian almost para-contact manifold}. In particular, if the metric $g$ is positive definite, then $(\epsilon)$-almost para-contact metric manifold is the usual {\it almost para-contact metric manifold} \cite{IS}.\\
The equation (\ref{2.2}) implies that 
\begin{equation}
\label{2.3}
g(X, \phi Y) = g(\phi X, Y)\,\,\ {\it and} \,\,\ g(X, \xi)=\epsilon \eta(X).
 \end{equation} 
 From (\ref{2.1}) and (\ref{2.3}) it follows that 
 \begin{equation}
\label{2.4}
g(\xi, \xi)=\epsilon.
 \end{equation} 
An $(\epsilon)$-almost para-contact metric structure is called an $(\epsilon)$-para-Sasakian structure if 
 \begin{equation}
\label{2.5}
(\nabla_{X}\phi)Y = -g(\phi X, \phi Y)\xi-\epsilon \eta(Y)\phi^2 X. \,\,\, X, Y \in TM,
 \end{equation} 
 where, $\nabla$ is the Levi-Civita connection with respect to $g$. A manifold endowed with $(\epsilon)$-para-Sasakian structure is called an {\it $(\epsilon)$-para-Sasakian manifold} \cite{TKYK}.
\par For $\epsilon=1$ and $g$ Riemannian, $M$ is the usual para-Sasakian manifold \cite{IS1}, \cite{SS1}. For $\epsilon=-1$, $g$ Lorentzian and $\xi$ replace by $-\xi$, $M$ becomes a Lorentzian para-Sasakian manifold \cite{KM}.
\par For an $(\epsilon)$-para-Sasakian manifold, it is easy to prove that  
\begin{eqnarray}\label{2.6}
 R(X,Y)\xi &=& \eta(X)Y-\eta(Y)X,\\
 \label{2.7}
 R(\xi,X)Y &=& \eta(Y)X-\epsilon g(X,Y)\xi, \\
\label{2.8}
 R(\xi,X)\xi &=& X-\eta(X)\xi, \\
\label{2.9}
 S(X, \xi) &=& -(n-1)\eta(X),\\
\label{2.10}
\nabla_{X}\xi &=& \epsilon \phi X.
\end{eqnarray}
For detail study of $(\epsilon)$-para-Sasakian manifold, see \cite{TKYK}. 
\par It is known that in a 3-dimensional $(\epsilon)$-para-Sasakian manifold (or, an $(\epsilon)$-para-Sasakian 3-manifold), the Riemannian curvature tensor and the Ricci tensor has the following form \cite{PKTK}:
\begin{eqnarray}\label{2.11}
R(X, Y)Z &=& \left(\frac{r}{2}+2\epsilon \right)\{g(Y, Z)X - g(X, Z)Y \}\nonumber\\
&-& \left(\frac{r}{2}+3\epsilon \right)\{g(Y, Z)\eta(X)\xi - g(X, Z)\eta(Y)\xi \nonumber\\
&+& \epsilon \eta(Y)\eta(Z)X - \epsilon \eta(X)\eta(Z)Y \},\\
\label{2.12} S(X, Y) &=& \left(\frac{r}{2} + \epsilon \right) g(X, Y) - \epsilon \left(\frac{r}{2} + 3 \epsilon \right) \eta(X)\eta(Y).
\end{eqnarray}
The projective curvature tensor $\mathcal{P}$ in a $(\epsilon)$-para-Sasakian 3-manifold $M$ is defined by
\begin{equation}\label{2.13}
\mathcal{P}(X, Y)U = R(X, Y)U - \frac{1}{2}[S(Y, U)X - S(X, U)Y].
\end{equation}
In an $(\epsilon)$-para-Sasakian 3-manifold $M$, the $\mathcal{Z}$ tensor takes the form
\begin{equation}\label{2.14}
\mathcal{Z}(X, Y) = \left(\frac{r}{2} + \epsilon + \psi \right) g(X, Y) - \epsilon \left(\frac{r}{2} + 3 \epsilon \right) \eta(X)\eta(Y).
\end{equation}
and scalar $\mathcal{Z}$ takes the form
\begin{equation}
\mathcal{Z} = \left(\frac{r}{2} + \epsilon + \psi \right)3 - \left(\frac{r}{2} + 3 \epsilon \right) = r + 3 \psi. \nonumber
\end{equation}
Also,
\begin{equation}\label{2.15}
\mathcal{Z}(X, \xi) = \left(\epsilon\psi - 2 \right)\eta(X).
\end{equation}
\section{\bf $\mathcal{Z}$-symmetries of $(\epsilon)$-para-Sasakian 3-manifolds}
In this section, we characterize, $\mathcal{Z}$-symmetric; $\mathcal{Z}$-semisymmetric; $\mathcal{Z}$-pseudosymmetric; and projectively $\mathcal{Z}$-semisymmetric conditions on an $(\epsilon)$-para-Sasakian 3-manifolds.
\\
We begin with the following:\\
\begin{defi}\label{def3.1}
A semi-Riemannian manifold $M$ is called locally symmetric if its curvature tensor $R$ is para-llel, that is, $\nabla R = 0$, where $\nabla$ denotes the Levi-Civita connection. As a proper generalization of locally symmetric manifolds, the notion of semisymmetric manifolds was defined by
\begin{equation}
R(X, Y)\cdot R = 0, \nonumber
\end{equation}
for any $X, Y \in TM$. A complete intrinsic classification of these spaces was given by Szabo \cite{ZIS}.
\end{defi}
\begin{defi}\label{def3.2}
A semi-Riemannian manifold $M$ is said to be $\mathcal{Z}$-symmetric if $\nabla\mathcal{Z} = 0$, and it is called $\mathcal{Z}$-semisymmetric if
\begin{equation}\label{3.1}
R(X, Y)\cdot\mathcal{Z} = 0,
\end{equation}
for any $X, Y\in TM$, where $R(X, Y)$ acts as a derivation on $\mathcal{Z}$.
\end{defi}
\par Let $M$ be a $\mathcal{Z}$-semisymmetric $(\epsilon)$-para-Sasakian 3-manifold. Then from (\ref{3.1}), we have
\begin{equation}
\mathcal{Z}(R(X, Y)U, V) + \mathcal{Z}(U, R(X, Y)V) = 0.\nonumber
\end{equation}
In particular, 
\begin{equation}\label{3.2}
\mathcal{Z}(R(\xi, Y)U, V) + \mathcal{Z}(U, R(\xi, Y)V) = 0.
\end{equation}
From (\ref{2.7}), we obtain
\begin{equation}\label{3.3}
\mathcal{Z}(R(\xi, Y)U, V) = -\epsilon g(Y, V)\mathcal{Z}(U, \xi) + \mathcal{Z}(U, Y)\eta(V)
\end{equation}
and
\begin{equation}\label{3.4}
 \mathcal{Z}(U, R(\xi, Y)V) = -\epsilon g(Y, U)\mathcal{Z}(V, \xi) + \mathcal{Z}(V, Y)\eta(U).
\end{equation}
The equations (\ref{3.2}), (\ref{3.3}) and (\ref{3.4}) together give
\begin{eqnarray}\label{3.5}
&& - \epsilon \{ g(Y, U)\mathcal{Z}(V, \xi) + g(Y, V)\mathcal{Z}(U, \xi)\} \nonumber\\
&& + \mathcal{Z}(V, Y)\eta(U) + \mathcal{Z}(U, Y)\eta(V) = 0.
\end{eqnarray}
Setting $V = \xi$ in (\ref{3.5}) and using (\ref{2.15}), we have
\begin{equation}\label{3.6}
\mathcal{Z}(U, Y) = (\psi-2\epsilon) g(Y, U).
\end{equation}
By making use of (\ref{2.14}) in (\ref{3.6}), we obtain
\begin{equation}\label{3.10}
\left(\frac{r}{2}+3\epsilon \right)[g(Y, U)-\epsilon \eta(Y)\eta(U)] = 0.
\end{equation} 
Since $g(Y, U)- \epsilon \eta(Y)\eta(U) = g(\phi Y, \phi U) \neq 0$, in general, therefore we obtain from (\ref{3.10}) that $\frac{r}{2}+3\epsilon = 0$, that is, 
\begin{equation}\label{3.11}
r = -6\epsilon.
\end{equation}
Next, using (\ref{3.11}) in (\ref{2.12}) we get
\begin{equation}\label{3.12}
S(X, Y) = -2\epsilon g(X, Y).
\end{equation}
That is, $M$ is an Einstein manifold. 
\\
\indent Conversely, suppose that the manifold $M$ be an Einstein. Then, from (\ref{1.1}) and (\ref{3.12}) we have (\ref{3.6}). Next, consider
\begin{eqnarray}\label{3.7}
R(X, Y)\cdot \mathcal{Z}(U, V) = \mathcal{Z}(R(X, Y)U, V) + \mathcal{Z}(U, R(X, Y)V).
\end{eqnarray}
By using (\ref{3.6}) in (\ref{3.7}) we obtain
\begin{eqnarray}\label{3.8}
R(X, Y)\cdot \mathcal{Z}(U, V) = (\psi - 2\epsilon)\{g(R(X, Y)U, V) + g(U, R(X, Y)V)\}.
\end{eqnarray}
It is known that, in an $(\epsilon)$-para-Sasakian manifold, the following relation holds:
\begin{equation}
\label{3.9} (R(X, Y)U, V) = - g(R(X, Y)V, U). 
\end{equation}
From (\ref{3.8}) and (\ref{3.9}), it follows that
\begin{equation}
R(X, Y)\cdot \mathcal{Z}(U, V) = 0.\nonumber
\end{equation}
That is, $M$ is $\mathcal{Z}$-semisymmetric. Hence, we are able to state the following result:
\begin{thm}\label{th3.1}
An $(\epsilon)$-para-Sasakian manifold is $\mathcal{Z}$-semisymmetric if and only if it is an Einstein manifold.
\end{thm}
Further, the manifold $M$ is an Einstein, implies trivially that $M$ is Ricci-symmetric.
\par Conversely, if $M$ is Ricci-symmetric, that is, $\nabla S = 0$.\\
In particular, 
\begin{equation}
(\nabla_X S)(Y, \xi) = \epsilon S(\phi X, Y) + 2 g(\phi X, Y) = 0.\nonumber
\end{equation} 
Replacing $X$ by $\phi X$ in the above equation, shows that the manifold is an Einstein manifold.
Therefore, by taking into account of theorem \ref{th3.1}, we state the following:
\begin{thm}\label{th3.2}
An $(\epsilon)$-para-Sasakian 3-manifold is $\mathcal{Z}$-semisymmetric if and only if it is Ricci-symmetric.
\end{thm}
Moreover, suppose that $M$ be Ricci-semisymmetric, that is, 
\begin{equation}
(R(X, Y)\cdot S)(U, V) = 0.\nonumber
\end{equation}
In particular,
\begin{equation}
(R(\xi, Y)\cdot S)(U, \xi) = 0,\nonumber
\end{equation}
this implies that
\begin{equation}
S(R(\xi, Y)U, \xi) + S(U, R(\xi, Y)\xi) = 0,\nonumber
\end{equation}
which in view of (\ref{2.7}) and (\ref{2.12}) gives (\ref{3.12}). 
\par Conversely, if $M$ is an Einstein manifold, then obviously, it is Ricci-semisymmetric. Thus, the manifold $M$ is Ricci-semisymmetric if and only if it is an Einstein manifold. Hence, by taking into account of theorem \ref{th3.1}, we have the following:
\begin{thm}\label{th3.3}
An $(\epsilon)$-para-Sasakian 3-manifold is $\mathcal{Z}$-semisymmetric if and only if it is Ricci-semisymmetric.
\end{thm}
\begin{cor}\label{cor3.1}
In an $(\epsilon)$-para-Sasakian 3-manifold, the following statements are equivalent:
\begin{enumerate}
    \item $M$ is an Einstein manifold.
    \item $M$ is Ricci-symmetric.
  \item $M$ is Ricci-semisymmetric.
  \item $M$ is $\mathcal{Z}$-semisymmetric.
  \end{enumerate} 
\end{cor}
It is clear that, $\nabla \mathcal{Z} = 0$ $\Rightarrow$ $ R\cdot \mathcal{Z} = 0$ $\Rightarrow$ $ \nabla S = 0$. Therefore, from corollary \ref{cor3.1}, we get:
\begin{cor}\label{cor3.2}
Every $\mathcal{Z}$-symmetric $(\epsilon)$-para-Sasakian 3-manifold is Ricci-symmetric.
\end{cor}
For a $(0, k)$-tensor field $T$ on $M$, $k\geq 1$, and a symmetric $(0, 2)$-tensor field $A$ on $M$, we define the $(0, k+2)$-tensor
fields $R\cdot T$ and $Q(A, T)$ by
\begin{eqnarray}
&& (R\cdot T)(X_1 ,...,X_k ; X, Y)\nonumber\\
&=&  -T(R(X, Y)X_1 , X_2 ,...,X_k ) -...-T(X_1 ,...,X_{k-1}, R(X,Y)X_k)\nonumber
\end{eqnarray} 
and
\begin{eqnarray}
&& Q(A, T)(X_1 ,...,X_k ; X, Y)\nonumber\\
&=& -T((X \wedge_A Y)X_1 , X_2 ,...,X_k ) -...-T(X_1 ,...,X_{k-1}, (X \wedge_A Y)X_k)\nonumber
\end{eqnarray}
respectively,
where $X\wedge_A Y$ is the endomorphism given by
\begin{equation}\label{3.15}
(X \wedge_A Y)Z = A(Y, Z)X - A(X, Z)Y.
\end{equation}
\begin{defi}\label{def3.3}
A semi-Riemannian $M$ is said to be pseudosymmetric (in the sense of R. Deszcz \cite{RD}) if
\begin{equation}
R\cdot R = L_{R}Q(g, R)\nonumber
\end{equation}
holds on $U_{R} = \{x \in M : R \neq 0 \frac{r}{n(n-1)}G\,\,\ at\,\,\ x\}$, where $G$ is the $(0, 4)$-tensor defined by $G(X_1, X_2, X_3, X_4) = g((X_{1}\wedge X_2)X_3 , X_4)$ and $L_{R}$ is some function on $U_{R}$.
\end{defi}
\begin{defi}\label{def3.4}
A semi-Riemannian manifold $M$ is said to be $\mathcal{Z}$-pseudosymmetric if 
\begin{equation}\label{3.16}
(R(X, Y)\cdot \mathcal{Z})(U, V) = L_{\mathcal{Z}}Q(g, \mathcal{Z})(U, V; X, Y)
\end{equation}
holds on the set $U_{\mathcal{Z}} = \{x \in M : \mathcal{Z} \neq 0\,\,\ at\,\,\ x\}$, and $L_{\mathcal{Z}}$ is some function on $U_{\mathcal{Z}}$. 
\end{defi}
\indent Let $M$ be a $\mathcal{Z}$-pseudosymmetric $(\epsilon)$-para-Sasakian 3-manifold. Then from (\ref{3.15}) and (\ref{3.16}), we have
\begin{equation}\label{3.17}
 (R(\xi, Y)\cdot \mathcal{Z})(U, \xi) = L_{\mathcal{Z}}\left[((\xi \wedge Y)\cdot \mathcal{Z})(U, \xi)\right].
\end{equation}
In an $(\epsilon)$-para-Sasakian 3-manifold, from (\ref{2.7}) and (\ref{3.15}) we get
\begin{equation}\label{3.18}
 R(\xi, X)Y = (-\epsilon)(\xi \wedge X)Y.
\end{equation}
In view of (\ref{3.17}) in (\ref{3.18}), it is easy to see that
\begin{equation}
 L_{\mathcal{Z}} = -\epsilon. \nonumber
\end{equation}
Hence, by taking into account of previous calculations and discussions, we conclude the following:
\begin{defi}\label{th3.4}
Let $M$ be an $(\epsilon)$-para-Sasakian 3-manifold. If $M$ is $\mathcal{Z}$-pseusosymmetric, then either $M$ is an Einstein manifold or $L_{\mathcal{Z}} = -\epsilon$ holds on $M$.
\end{defi}
\indent If $L_{\mathcal{Z}} \neq -\epsilon$, then immediately, we obtain the following:
\begin{cor}\label{cor3.3}
Every $\mathcal{Z}$-pseudosymmetric $(\epsilon)$-para-Sasakian 3-manifold with $L_{\mathcal{Z}} \neq -\epsilon$ is an Einstein manifold.
\end{cor}
\indent But $L_{\mathcal{Z}}$ need not be zero, in general and hence there exists $\mathcal{Z}$-pseudosymmetric manifolds which are not $\mathcal{Z}$-semisymmetric. Thus the class of $\mathcal{Z}$-pseudosymmetric manifolds is a natural extension of the class of $\mathcal{Z}$-semisymmetric manifolds. Thus, if $L_{\mathcal{Z}} \neq 0$ then it is easy to see that $R\cdot \mathcal{Z} = (-\epsilon) Q(g, \mathcal{Z})$, which implies that the pseudosymmetric function $L_{\mathcal{Z}} = -\epsilon$. Therefore, we able to state the following result:
\begin{cor}\label{cor3.4}
 Every $(\epsilon)$-para-Sasakian 3-manifold is $\mathcal{Z}$-pseudosymmetric of the form $R\cdot \mathcal{Z} = (-\epsilon) Q(g, \mathcal{Z})$.
 \end{cor}
\begin{defi}\label{def3.5}
A semi-Riemannian manifold $M$ is said to be projectively $\mathcal{Z}$-semisymmetric if 
\begin{equation}\label{3.20}
\mathcal{P}(X, Y)\cdot\mathcal{Z} = 0,
\end{equation}
for any $X, Y \in TM$, where $\mathcal{P}$ is the projective curvature tensor. 
\end{defi}
\indent Let $M$ be a projectively $\mathcal{Z}$-semisymmetric $(\epsilon)$-para-Sasakian 3-manifold. Then from (\ref{3.20}), we have
\begin{equation}
\mathcal{Z}(\mathcal{P}(X, Y)U, V) + \mathcal{Z}(U, \mathcal{P}(X, Y)V) = 0.\nonumber
\end{equation}
This implies
\begin{equation}\label{3.21}
\mathcal{Z}(\mathcal{P}(\xi, Y)U, V) + \mathcal{Z}(U, \mathcal{P}(\xi, Y)V) = 0.
\end{equation}
In an $(\epsilon)$-para-Sasakian 3-manifold, from (\ref{2.13}) we have 
\begin{equation}\label{3.21a}
\mathcal{P}(\xi, X)Y = -\frac{1}{2}S(X, Y)\xi - \epsilon g(X, Y)\xi.
\end{equation}
Using (\ref{2.7}) and (\ref{3.21a}), we obtain
\begin{equation}\label{3.22}
\mathcal{Z}(\mathcal{P}(\xi, Y)U, V) = \left[-\frac{1}{2}S(Y, U) - \epsilon g(Y, U)\right]\mathcal{Z}(V, \xi)
\end{equation}
and
\begin{equation}\label{3.23}
 \mathcal{Z}(U, \mathcal{P}(\xi, Y)V) = \left[-\frac{1}{2}S(Y, V) - \epsilon g(Y, V)\right]\mathcal{Z}(U, \xi).
\end{equation}
The equations (\ref{3.21}), (\ref{3.22}) and (\ref{3.23}) together give
\begin{eqnarray}\label{3.24}
&& \left[\frac{1}{2}S(Y, U) + \epsilon g(Y, U)\right]\mathcal{Z}(V, \xi) + \left[\frac{1}{2}S(Y, V) + \epsilon g(Y, V)\right]\mathcal{Z}(U, \xi) = 0.
\end{eqnarray}
Setting $V = \xi$ in (\ref{3.24}) and then using (\ref{2.12}) and (\ref{2.15}), we have
\begin{equation}\label{3.25}
(\psi - 2\epsilon)[S(U, Y) + 2\epsilon g(Y, U)]=0.
\end{equation}
This implies either $\psi = 2\epsilon$ or $S(U, Y) = -2\epsilon g(Y, U)$.
\indent If $\psi = 2\epsilon$, then from (\ref{2.15}) we obtain (\ref{3.10}). It shows that $M$ is an Einstein manifold. 
\indent Therefore, in both of the cases, manifold $M$ reduces to an Einstein manifold. Hence, we state the following:
\begin{thm}\label{th3.5}
Every projectively $\mathcal{Z}$-semisymmetric $(\epsilon)$-para-Sasakian 3-manifold is an Einstein manifold.
\end{thm}

\bibliographystyle{spphys}       

\end{document}